\documentclass[graybox]{svmult}

\usepackage{mathptmx}
\usepackage{helvet}
\usepackage{courier}
\usepackage{type1cm}
\usepackage{makeidx}
\usepackage{graphicx}

\usepackage{multicol}
\usepackage[bottom]{footmisc}

\usepackage{tikz}

\usepackage{url}

\usepackage{amsfonts}
\usepackage{amsmath}
\usepackage{amssymb}

\begin{document}

\title{On parameter-dependent inhomogeneous boundary-value problems in Sobolev spaces}
\titlerunning{On parameter-dependent boundary-value problems}
\author{Olena Atlasiuk, Vladimir Mikhailets, and Jari Taskinen}
\institute{Olena Atlasiuk \at University of Helsinki, Department of Mathematics and Statistics, P.O. Box 68, Pietari Kalmin katu 5, 00014 Helsinki, Finland and Institute of Mathematics of the National Academy of Sciences of Ukraine, st. Tereschenkivska 3, 01024 Kyiv, Ukraine \at \email{olena.atlasiuk@helsinki.fi}
\and Vladimir Mikhailets \at King's College London, Strand, WC2R 2LS London, UK and Institute of Mathematics of the National Academy of Sciences of Ukraine, st. Tereschenkivska 3, 01024 Kyiv, Ukraine
\and Jari Taskinen \at University of Helsinki, Department of Mathematics and Statistics, P.O. Box 68, Pietari Kalmin katu 5, 00014 Helsinki, Finland }
\maketitle

\abstract{
We study a wide class of linear inhomogeneous boundary-value problems for $r$th order ODE-systems depending on a
parameter $\mu$ belonging to a general metric space $\mathcal M$. The solutions belong to
the Sobolev spaces $(W^{n+r}_p)^m$, $n\in\mathbb{N}\cup\{0\}$, $m, r \in \mathbb{N}$, $1\leq p\leq \infty$.
The boundary conditions are of a most general form $By=c$, where $B$ is an arbitrary continuous operator
from $(W^{n+r}_p)^m$ to $\mathbb{C}^{rm}$. Thus, they  may contain derivatives of the unknown vector
function of integer and/or fractional orders $\geq r$. We find necessary and
sufficient conditions for the continuity of solutions with respect to the parameter $\mu$. We also prove
that the solutions of the original problems can be approximated in the space $(W^{n+r}_p)^m$ by
solutions of ODE-systems with polynomial coefficients, right-hand sides of the equation, and multipoint boundary conditions,
which are independent of the original problem's right-hand sides.
}


\section{Introduction}
\label{sec:1}

In some applied problems, boundary-value problems naturally arise in which the boundary conditions contain derivatives whose order exceeds the order of the differential equation. The theory of such problems is not yet sufficiently developed and is concentrated mainly around elliptic boundary-value problems in Sobolev spaces. Such problems are Fredholm in the corresponding pairs of Sobolev spaces, and the question of their index is quite complicated and is studied in Boutet de Monvel’s theory. The simplest class of elliptic differential equations is a linear ODE.

It is well known that for systems of linear ODEs, the solution of the Cauchy problem always exists and is unique. Unlike them, solutions of boundary-value problems may not exist even under the simplest two-point boundary conditions. Until recently, the widest  class of inhomogeneous boundary-value problems was the class of problems with general boundary conditions of the form
$Ly=f$,  $By=c$,
where $L$ is a system of differential equations of order $r
\in \mathbb{N}$ with Lebesgue summable coefficients, and $B$ is a linear continuous finite-dimensional operator $B \colon W^r_1 \rightarrow \mathbb{C}^r$. Therefore, the general boundary conditions may contain derivatives of an unknown function of order $\leq r-1$.

In this paper, we describe and investigate a broader class of boundary-value problems for differential systems of arbitrary order on a finite interval. The coefficients and right-hand sides of these problems belong to some Sobolev space, due to which their solutions have additional smoothness in the Sobolev scale. This allows us to correctly define for them the widest possible class of inhomogeneous boundary conditions, which may contain derivatives of an unknown function of integer and/or fractional order, which may exceed the order of the differential equation.

The main results of the article are as follows:
\begin{itemize}
\item[---] Indices and d-characteristics of the introduced boundary-value problems are found.
\item[---] Constructive necessary and sufficient conditions for the continuity of solutions in a parameter from an abstract metric space are obtained. This allows us to cover the cases of discrete, continuous, and functional parameters from a single point of view.
\item[---] It is proved that solutions continuous in a parameter are stable in the sense of Ulam--Hyers.
\item[---] It is established that the solutions of an arbitrary solvable boundary-value problem are the limits of solutions of problems with polynomial coefficients and multipoint boundary conditions in the corresponding Sobolev spaces.
\end{itemize}

\section{General results on the continuous parameter dependence}
\label{sec:2}

Let $(a,b)\subset\mathbb{R}$ be a finite interval and let the parameters $
n\in\mathbb{N}\cup\{0\}$, $\{m, r, k\}\subset\mathbb{N}$, $1\leq p\leq \infty$ be arbitrary. We denote by
\begin{gather*}
W_p^{n+r}\bigl([a,b];\mathbb{C})
:= \bigl\{y\in C^{n+r-1}([a,b];\mathbb{C})\colon y^{(n+r-1)}\in AC[a,b], \, y^{(n+r)}\in L_p[a,b]\bigr\}
\end{gather*}
the usual  complex Sobolev space and we set $W_p^0:=L_p$. This space is Banach with respect to the norm
$$
\bigl\|y\bigr\|_{n+r,p}=\sum_{s=0}^{n+r}\bigl\|y^{(s)}\bigr\|_{p},
$$
where $\|\cdot\|_p$ stands for the norm in the Lebesgue space $L_p\bigl([a,b]; \mathbb{C}\bigr)$. We
use the Sobolev spaces
$
(W_p^{n+r})^{m}:=W_p^{n+r}\bigl([a,b];\mathbb{C}^{m}\bigr)$ and $(W_p^{n+r})^{m\times m}:=W_p^{n+r}\bigl([a,b];\mathbb{C}^{m\times m}\bigr),$
which consist, respectively, of vector- and matrix-valued functions with ele\-ments belonging to $W_p^{n+r}$.
The norms in these spaces are defined to be the sums of the Sobolev-norms of the components, and the
same notation $\|\cdot\|_{n+r,p}$ is used in all cases, which will be clear from the context.
The same convention will be applied to all other Banach spaces. 

We consider a linear boundary-value problem of the form
\begin{equation}\label{bound_pr_1}
(Ly)(t):=y^{(r)}(t) + \sum\limits_{\ell=1}^rA_{r-\ell}(t)y^{(r-\ell)}(t)=f(t), \quad t\in(a,b),\\
\end{equation}
\begin{equation}\label{bound_pr_2}
By=c.
\end{equation}
We suppose that matrix-valued functions $A_{r-\ell}(\cdot) \in (W_p^n)^{m\times m}$, a vector-valued function $f(\cdot) \in (W^n_p)^m$, a vector $c \in \mathbb{C}^{rm}$, a linear continuous operator
\begin{equation}\label{oper_B_class}
B\colon (W^{n+r}_p)^m\rightarrow\mathbb{C}^{rm}
\end{equation}
are arbitrarily chosen and that the vector function $y(\cdot)\in (W_{p}^{n+r})^m$ is unknown.

The boundary condition \eqref{bound_pr_2} consists of $rm$ scalar conditions for a system of $m$ differential equations of $r$-th order, where we represent vectors and vector-valued functions as columns. A solution to the boundary-value problem \eqref{bound_pr_1}, \eqref{bound_pr_2} is understood as a vector-valued function $y(\cdot)\in (W_{p}^{n+r})^m$ that satisfies both equation \eqref{bound_pr_1} (everywhere if $n\geq 1$, and almost everywhere if $n=0$) on $(a,b)$ and equality \eqref{bound_pr_2}. If the parameter $n$ increases, so does the class of linear operators \eqref{oper_B_class}. When $n=0$, this class contains all operators that set the general boundary conditions.

The set of solutions to equation \eqref{bound_pr_1} coincides with the space $(W_{p}^{n+r})^m$ when the right-hand side of the equation runs over $(W_{p}^{n})^m$. Hence, the boundary condition \eqref{bound_pr_2} with the operator of the form \eqref{oper_B_class} is the most general for this equation.

We rewrite the inhomogeneous boundary-value problem \eqref{bound_pr_1}, \eqref{bound_pr_2} in the form of a linear operator equation
$
(L,B)y=(f,c).
$
Here, $(L,B)$ is a bounded linear operator on the pair of Banach spaces
\begin{equation}\label{(L,B)}
(L,B)\colon (W^{n+r}_p)^m\rightarrow (W^{n}_p)^m \oplus \mathbb{C}^{rm},
\end{equation}
which follows from the definition of the Sobolev spaces involved and from the fact that $W_p^n$ is a Banach algebra.

\begin{theorem}\label{opBinSobolev}
Let $1\leq p \leq\infty$ and $1/p + 1/p'=1$, and let  $t_0 \in [a, b]$, the matrix $\big(\alpha_{s}\big)_{s=1}^{n+1-r} \subset \mathbb{C}^{rm\times rm}$ and the matrix-valued function $\Phi(\cdot)\in L_{p'}([a, b] ; \mathbb{C}^{rm\times rm})$ be given.
\begin{itemize}
\item [(1)] \, The operator $B$ defined by
\begin{equation} \label{diyaBnash}
By=\sum _{s=0}^{n+r-1} \alpha_{s}\,y^{(s)}(t_0)+\int_{a}^b \Phi(t)y^{(n+r)}(t){\rm d}t, \quad y(\cdot)\in (W_{p}^{n+r})^{m},
\end{equation}
acts continuously from $(W_{p}^{n+r})^{m}$ into $\mathbb{C}^{rm}$. 
\item [(2)] \, If $p \neq \infty$, then every bounded operator $B\colon (W^{n+r}_p)^m\rightarrow\mathbb{C}^{rm}$ admits a unique canonical representation of the form \eqref{diyaBnash}.
\end{itemize}
\end{theorem}

It should be noted that in the case of $p =\infty$, not all operators $B$ can be presented in the form
\eqref{diyaBnash}, since there are continuous operators $B$ that are defined by integrals over finitely additive
measures (see, for instance, \cite{Bhaskara,Dunford,KantAk1982}).

Let $E_{1}$ and $E_{2}$ be Banach spaces. A linear bounded operator $T\colon E_{1}\rightarrow E_{2}$ is called a Fredholm operator if its kernel and co-kernel are finite-dimen\-si\-onal. If $T$ is a Fredholm operator, then its range $T(E_{1})$ is closed in $E_{2}$, and its index is finite
$
\mathrm{ind}\,T:=\dim\ker T-\dim\big(E_{2}/T(E_{1})\big)\in \mathbb{Z} = \{0, \pm 1, \pm2 , \ldots \}$ (see, e.g., \cite[Lemma~19.1.1]{Hermander1985}).

\begin{theorem}\label{th_fredh high}
The bounded linear operator \eqref{(L,B)} is a Fredholm one with zero index.
\end{theorem}

Theorem \ref{th_fredh high} naturally raises the question of finding d-characteristics of the operator $\big(L,B\big)$, that is,  $\operatorname{dim} \operatorname{ker}\big(L,B\big)$ and $\operatorname{dim} \operatorname{coker}\big(L,B\big)$. This is a quite difficult task because d-characteristics may vary even with arbitrarily small one-dimensional additive perturbations.

To formulate the following result, let us introduce some notation and definitions.

For each number $i \in \{1,\dots, r\}$, we consider the family of matrix Cauchy problems with the initial conditions:
\begin{equation}\label{zad kosh1}
Y_i^{(r)}(t)+\sum\limits_{j=1}^rA_{r-j}(t)Y_i^{(r-j)}(t)=O_{m},\quad t\in (a,b),
\end{equation}
\begin{equation}\label{zad kosh2}
Y_i^{(j-1)}(a) = \delta_{i,j}I_m,\quad j \in \{1,\dots, r\},
\end{equation}
where $Y_i(\cdot)$ is an unknown $(m\times m)$--matrix-valued function. As usual, $O_{m}$ stands for the zero $(m\times m)$--matrix, $I_{m}$ denotes the identity $(m\times m)$--matrix, and $\delta_{i,j}$ is the Kronecker delta. Each Cauchy problem \eqref{zad kosh1}, \eqref{zad kosh2} has a unique solution $Y_i\in(W_p^{n+r})^{m\times m}$ due to \cite[Lemma 4.1]{NovAM2023}. Certainly, if $r=1$, we use the designation  $Y(\cdot)$ for $Y_1(\cdot)$.

Let $\left[BY_i\right]$ denote the number $(m\times m)$--matrix whose $j$-th column is the result of the action of $B$ on the $j$-th column of the matrix-valued function~$Y_i$.

\begin{definition}\label{defin_harm}
A bloc square number matrix
\begin{equation}\label{matrix_BY}
M(L,B):=\big(\left[BY_1\right],\dots,\left[BY_{r}\right]\big) \in \mathbb{C}^{rm\times rm}
\end{equation}
is called the characteristic matrix of the inhomogeneous boundary-value problem \eqref{bound_pr_1}, \eqref{bound_pr_2}. Note that this matrix consists of $r$ square block columns $\left[BY_k\right]\in \mathbb{C}^{m\times m}$.
\end{definition}


\begin{theorem}\label{th dimker}
The dimensions of the kernel and co-kernel of the operator \eqref{(L,B)} are equal to the dimensions of the kernel and co-kernel of the characteristic matrix \eqref{matrix_BY}, respectively.
\end{theorem}

Theorem \ref{th dimker} implies the following necessary and sufficient conditions for the invertibility of \eqref{(L,B)}:

\begin{corollary}\label{invertible}
The operator \eqref{(L,B)} is invertible if and only if the square matrix $M(L,B)$ is nonsingular \cite{NovAM2023}.
\end{corollary}


Let $\mathcal{M}$ be an arbitrary metric space. In the sequel, $\mu \in \mathcal{M}$ will denote a free
parameter whereas $\mu_0 \in \mathcal{M}$ denotes an arbitrarily  fixed one.
We consider the following linear boundary-value problem of the form \eqref{bound_pr_1}, \eqref{bound_pr_2} for an unknown
vector-valued function $y(\cdot,\mu)\in (W_{p}^{n+ r})^m$,
\begin{equation}\label{bound_z1}
\left(L(\mu)y(\mu)\right)(t):=y^{(r)}(t,\mu) + \sum\limits_{\ell=1}^rA_{r-\ell}(t,\mu)y^{(r-\ell)}(t,\mu)=f(t,\mu),
\end{equation}
\begin{equation} \label{bound_z2}
B(\mu)y(\mu)=c(\mu), \quad t\in(a,b),
\end{equation}
where matrix-valued functions $A_{r-\ell}(\cdot,\mu) \in (W_p^n)^{m\times m}$, a vector-valued function
$f(\cdot,\mu)$ $\in (W^n_p) ^m$, a vector $c(\mu) \in \mathbb{C}^{rm}$ and a
continuous linear operator
\begin{equation}\label{oper_B(e)v}
B(\mu)\colon (W^{n+r}_p)^m\rightarrow\mathbb{C}^{rm}
\end{equation}
are given and arbitrary.
The boundary condition \eqref{bound_z2} consists of $rm$ scalar conditions for a system of $m$ differential
equations of $r$-th order. 

A solution to problem \eqref{bound_z1}, \eqref{bound_z2} is understood as a vector-valued function
$y \in (W_{p}^{n+r})^m$ which satisfies both equation \eqref{bound_z1} (everywhere if $n\geq 1$, and
almost everywhere if $n=0$) on $(a,b)$ and equality \eqref{bound_z2}.
As explained above, we refer to the general boundary condition \eqref{bound_z2} with an arbitrary continuous
operator \eqref{oper_B(e)v} as  \textit{generic} for the differential system \eqref{bound_z1}. It covers
all classical types of boundary conditions, such as initial conditions in the Cauchy problem, various
multipoint conditions, integral conditions,  mixed boundary conditions, as well as non-classical conditions
containing fractional derivatives, where the order of the derivatives may exceed the order of the
differential equation. Finally, we do not pose any {\it a priori} assumption on the regularity of
the matrix-value functions $A_{r-\ell}(t,\mu)$ with respect to $\mu$.

We write problems \eqref{bound_z1}, \eqref{bound_z2} in the form of a linear operator equation
$$ \big(L(\mu),B(\mu)\big)y(\mu)=\big(f(\mu),c(\mu)\big), $$
where $\big(L(\mu),B(\mu)\big)$ is the family of continuous linear operators
\begin{equation}\label{(L,B)vp}
\big(L(\mu),B(\mu)\big) \colon (W^{n+r}_p)^m\to (W^{n}_p)^m\times\mathbb{C} ^{rm}.
\end{equation}

According to Theorem \ref{th_fredh high}, all operators in the family  \eqref{(L,B)vp} are Fredholm with index zero for every  $\mu$.

\begin{definition}\label{defin_vp}
We say that a solution to the boundary-value problem \eqref{bound_z1}, \eqref{bound_z2} depends continuously
on the parameter $\mu$ at a limit point $\mu_0$ of the metric space $\mathcal{M}$,
if the following two conditions are satisfied:
\begin{itemize}
\item[$(\ast)$] \,\, There exists a positive number $\varepsilon$ such that, for all $\mu\in \mathcal{B}(\mu_0, \varepsilon)$ and arbitrary right-hand sides $f(\cdot;\mu)\in (W^{n}_p)^{m}$ and $c(\mu)\in\mathbb{C}^{rm}$, the problem has a unique solution $y(\cdot;\mu)$ in the space  $(W^{n+r}_p)^{m}$;
\item [$(\ast\ast)$] \,\, The convergence of the right-hand sides $f(\cdot;\mu)\to f(\cdot;\mu_0)$ in $(W_p^n)^{m}$ and $c(\mu)\to c(\mu_0)$ in $\mathbb{C}^{rm}$ as $\mu\to\mu_0$ implies the convergence of the solutions
\begin{equation*}\label{4.guv}
y(\cdot,\mu)\to y(\cdot,\mu_0)\quad\mbox{in}\quad (W^{n+r}_p)^{m} \quad\mbox{as}\quad\mu\to\mu_0.
\end{equation*}
\end{itemize}
\end{definition}

Throughout this article, we will assume that the following condition (0) for the point
$\mu_0 \in \mathcal{M}$ is fulfilled.

\medskip

\noindent \textbf{Condition (0)}. {\it The homogeneous boundary-value problem has only a trivial solution
$L(\mu_0)y(t,\mu_0)=0$, $t\in (a,b)$, $B(\mu_0)y(\cdot,\mu_0)=0$}.

\medskip

\noindent We will also consider the following two conditions on the left-hand sides of the problem
\eqref{bound_z1}, \eqref{bound_z2}.

\medskip

\noindent \textbf{Limit Conditions} as $\mu\to\mu_0$:

\smallskip

\noindent
(I) {\it $A_{r-\ell}(\cdot;\mu)\to A_{r-\ell}(\cdot;\mu_0)$ in the space $(W^{n}_p)^{m\times m}$ for every  $\ell\in\{1,\ldots, r\}$;}

\smallskip

\noindent
(II) {\it $B(\mu)y\to B(\mu_0)y$ in the space $\mathbb{C}^{m}$ for all $y\in(W^{n+r}_p)^m$.
}

\medskip
Now, we can formulate necessary and sufficient conditions for the continuity of the solutions to
the boundary-value problem \eqref{bound_z1}, \eqref{bound_z2} with respect to an abstract parameter.

\begin{theorem}\label{nep v}
The solution to the boundary-value problem \eqref{bound_z1}, \eqref{bound_z2} depends continuously on the parameter~$\mu$ at $\mu_0\in \mathcal{M}$ if and only if this problem satisfies Condition \textup{(0)} and Limit Conditions~\textup{(I)} and~\textup{(II)}.
\end{theorem}

\begin{corollary}
If Condition \textup{(0)} and Limit Conditions~\textup{(I)} and~\textup{(II)} are satisfied for all
$\mu \in \mathcal{M}$  and the right-hand sides $f$ and  $c$ are fixed,
then the solution to the boundary-value problem \eqref{bound_z1}, \eqref{bound_z2} exists and
is unique for every $\mu \in \mathcal{M}$ and  belongs to the space $C\big(\mathcal{M}; (W^{n+r}_p)^{m}\big)$.
\end{corollary}

It is worth noting that using an arbitrary metric space $\mathcal{M}$ in Theorem \ref{nep v} yields
a unified approach to both continuous and discrete parameters.

In the case of $r=1$, $\mathcal{M} = [0, \varepsilon_0]$, $\varepsilon_0>0$, $\mu_0=0$, Theorem \ref{nep v} was proved in  \cite[Theorem 1]{AtlasiukMikhailets20192} and in the case of $r=1$, $\mathcal{M} = I \subset \mathbb{R}$, where $I$ is an interval on $\mathbb{R}$, in  \cite[Theorem 1]{AtlasiukMikhailets2025}.

We supplement the previous result with an error estimate  for   the solution $y(\cdot;\mu)$,
namely, considering $y(\cdot;\mu)$ as an approximate solution of the problem \eqref{bound_z1} for the
parameter value $\mu_0$, \eqref{bound_z2}, we show that
$
\bigl\|y(\cdot;\mu_0)-y(\cdot;\mu)\bigr\|_{n+r,p}
$
is proportional to the discrepancy
$\widetilde{d}_{n,p}(\mu):=
\bigl\|L(\mu)y(\cdot;\mu_0)-f(\cdot;\mu)\bigr\|_{n,p}+
\bigl\|B(\mu)y(\cdot;\mu_0)-c(\mu)\bigr\|_{\mathbb{C}^{rm}}$.

\begin{theorem}\label{3.6.th-bound v}
Let $\mu_0 \in \mathcal{M}$ and assume that the boundary-value problem \eqref{bound_z1},
\eqref{bound_z2} satisfies Condition \textup{(0)} and Limit Conditions \textup{(I)} and \textup{(II)}.
Then there exist positive numbers $\varepsilon$, $\gamma_{1}$, and $\gamma_{2}$ such that
\begin{equation*}\label{3.6.bound}
\begin{aligned}
\gamma_{1}\,\widetilde{d}_{n,p}(\mu)
\leq\bigl\|y(\cdot;\mu_0)-y(\cdot;\mu)\bigr\|_{n+r,p}\leq
\gamma_{2}\,\widetilde{d}_{n,p}(\mu),
\end{aligned}
\end{equation*}
for all $\mu\in\mathcal{B}(\mu_0, \varepsilon)$. Here, $\varepsilon$, $\gamma_{1}$, and $\gamma_{2}$ do not depend on $y(\cdot;\mu_0)$ or  $y(\cdot;\mu)$.
\end{theorem}

The right-hand inequality in the formula \eqref{3.6.bound} can be interpreted as the Ulam--Hyers  stability of solutions of the boundary-value problem with parameter as  $\mu \rightarrow \mu_0$ (see, for example, \cite{Hyers,Rus}).

\section{Limit theorems concerning the operators of the boundary-value problem}
\label{sec:3}
Our aim in this section is to present conditions  on the coefficients of differential expressions and
operators $B(\mu)$ under which $\big(L(\mu),B(\mu)\big)$ converges to the operator $\big(L(\mu_0),B(\mu_0)\big)$
in the strong and uniform operator topologies.
First, we formulate necessary and sufficient conditions for the strong and uniform convergence of the family
of operators $L(\mu)$ to the operator $L(\mu_0)$.

\begin{theorem}\label{eqival ym L}
Let $1\leq p\leq \infty$. The following convergence conditions are equivalent, when $\mu\to\mu_0$
in the metric space $\mathcal{M}$:
\begin{itemize}
\item [(I)] \,\,\, $A_{r-\ell}(\cdot,\mu) \rightarrow A_{r-\ell}(\cdot, \mu_0)$ in the Banach space
$(W^{n}_p)^{m\times m}$ for all $\ell\in\{1,\ldots, r\}$;
\item [(II)] \,\, $L(\mu) \rightarrow L(\mu_0)$ in the uniform operator topology;
\item [(III)] \,\, $L(\mu) \rightarrow L(\mu_0)$ in the strong operator topology.
\end{itemize}
\end{theorem}

Note that in the case where the metric parameter $\mu$ is a natural number, Theorem \ref{eqival ym L} was
proved in \cite[Lemma 6.1]{NovAM2023}. The proof of the general case was
given in \cite[Theorem 3.1]{AtlMikhJ}.


We next formulate necessary and sufficient conditions for the strong and uniform convergence of the family  $B(\mu)$. To this end, we consider the following asymptotic conditions as $\mu\to\mu_0$.

\begin{enumerate}
\renewcommand{\labelenumi}{\alph{enumi})}
\renewcommand{\theenumi}{\alph{enumi})}
\item[(a)] $\alpha_s(\mu)\rightarrow\alpha_s(\mu_0)$ in $\mathbb{C}^{rm\times rm}$ for every
$s\in\{0,\dots, n+r-1\}$;
\item[(b)]
$\left\|\Phi(\cdot,\mu)\right\|_{p'}=O(1)$;
\item[(c)] $\int\limits_a^t \Phi(\tau,\mu)d\tau\rightarrow\int\limits_a^t\Phi(\tau, \mu_0)d\tau$
in the space $\mathbb{C}^{rm\times rm}$ for all $t\in(a,b]$;
\item[(d)] $\|\Phi(\cdot,\mu)-\Phi(\cdot, \mu_0)\|_{p'}\rightarrow0$.
\end{enumerate}

It is easy to see that condition $(d)$ is stronger than conditions $(b)$ and $(c)$.

\begin{theorem}\label{strongunifB}
Let $1\leq p<\infty$. The operators $B(\mu)$ converge strongly to the operator $B(\mu_0)$,
as  $\mu\to\mu_0$, if and only if conditions $(a)$, $(b)$ and $(c)$ hold. The convergence is
uniform if and only if $(a)$ and $(d)$ are satisfied.
\end{theorem}

\section{Approximation by solutions of multipoint boundary-value problems}
\label{sec:4}

Let us next apply the above results to the approximation of solutions of an inhomogeneous
boundary-value
problem by solutions of a sequence of boundary-value problems with polynomial coefficients, right-hand sides, and multipoint
boundary conditions. We will consider the case $p < \infty $, since the case $p=\infty$ is
significantly different and will be omitted here.

To  formulate the statement of the problem, we consider a well-posed boundary-value problem with inhomogeneous boundary conditions
\begin{equation}\label{bound_pr_aaaaa}
(L_0y_0)(t):=y^{(r)}_0(t) + \sum\limits_{\ell=1}^rA_{r-\ell,0}(t)y^{(r-\ell)}_0(t)=f_0(t), \quad t\in(a,b),
\end{equation}
\begin{equation}\label{bound_pr_bbbb}
B_0y_0:=\sum _{s=0}^{n+r-1} \alpha_{s,0}\,y^{(s)}_0(t_0)+\int_{a}^b \Phi_0(t)y^{(n+r)}_0(t){\rm d}t=c_0,
\end{equation}
when the matrix-valued functions $A_{r-\ell}(\cdot) \in (W_p^n)^{m\times m}$, the vector-valued function $f_0(\cdot) \in (W^n_p)^m$, the vector $c_0 \in \mathbb{C}^{rm}$, the numerical matrices $\alpha_{s,0} \in \mathbb{C}^{rm\times rm}$ and the matrix-valued function $\Phi_0(\cdot)\in L_{p'}([a, b] ; \mathbb{C}^{rm\times rm})$, $p \in [1, \infty)$, $p^{-1}+p'^{-1}=1$, are given. As we proved
in Theorem \ref{opBinSobolev}, an arbitrary inhomogeneous boundary condition for equation \eqref{bound_pr_aaaaa} admits a unique canonical representation of the form \eqref{bound_pr_bbbb}, where $t_0$ is an arbitrary fixed point of the interval $[a,b]$.

Consider simultaneously a sequence of multipoint boundary-value problems
\begin{equation}\label{bound_pr_aaaww}
(L_ky_k)(t):=y^{(r)}_k(t) + \sum\limits_{\ell=1}^rA_{r-\ell,k}(t)y^{(r-\ell)}_k(t)=f_k(t), \quad t\in(a,b),
\end{equation}
\begin{equation}\label{bound_pr_bbbww}
B_ky_k:=\sum _{s=0}^{n+r-1} \alpha_{s,k}\,y^{(s)}_k(t_0)+\sum\limits_{j=0}^{N(k)}{\beta_{j,k}
y^{(n+r-1)}_k(t_{j,k})}=c_k,  \ \ \ k=1,2,3, \ldots.
\end{equation}
Here, the elements of the matrix-valued functions $A_{r-\ell,k}$ belong to
some dense set  $\mathcal{F}$ in the space $(W_p^n)^{m\times m}$, $G$ is a dense set in the space $(W_p^n)^{m}$, and $f_k \in G$ for all $k$. Moreover, for all indices, $\alpha_{s,k}, \beta_{j,k} \in
\mathbb{C}^{m\times m}$ and the points $t_{j,k}$ belong to some dense
set $\mathcal{P}$ in the interval $[a,b]$, and we have
$
f_k \rightarrow f_0$, $c_k \rightarrow c_0$ as $k \rightarrow\infty.
$

We next study the natural problem of the existence of a sequence of boundary-value problems \eqref{bound_pr_aaaww}, \eqref{bound_pr_bbbww}, whose solutions satisfy the asymptotic formula
$$
y_k \rightarrow y_0, \quad \mbox{in}\quad (W_p^{n+r})^{m}  \ \ \mbox{as} \ k\rightarrow\infty.
$$
We will give a positive answer to this question, which is based on  Theorem \ref{nep v}
on the continuity of the solutions of boundary-value problems with respect to a parameter
belonging to an abstract metric space $\mathcal{M}$.

In fact, we set $\mathcal{M}=\mathbb{Z}_+$ and introduce a metric on $\mathcal{M}$ in such a way that 0 is the only limit point in the metric space $\big(\mathbb{Z}_+, d\big)$ and
$
d(0,n)\rightarrow 0 \Leftrightarrow n\rightarrow\infty
$ holds.

We present the main result of this section.

\begin{theorem}\label{thapro}
Assume that the homogeneous boundary-value problem has only a trivial solution. Then,
there exists a sequence of well-posed boundary-value problems of the form \eqref{bound_pr_aaaww}, \eqref{bound_pr_bbbww} with polynomial coefficients and right-hand sides such that
\begin{itemize}
\item[(i)] \, if $1<p<\infty$, then
$\big\| \left(L_k, B_k\right)^{-1}-\left(L_0, B_0\right)^{-1}\big\|\rightarrow 0$, as $k\rightarrow \infty$;
\item[(ii)] \, if $p=1$, then
$\left(L_k, B_k\right)^{-1}  \stackrel{s}{\longrightarrow} \left(L_0, B_0\right)^{-1}$ and $y_k \rightarrow y_0$ in  $(W_1^{n+r})^m$, as $k\rightarrow \infty$;
\item[(iii)] \,\, if $p=\infty$ and the boundary conditions admit the representation of the form \eqref{diyaBnash}, then
$\big\| \left(L_k, B_k\right)^{-1}-\left(L_0, B_0\right)^{-1}\big\|\rightarrow 0$, as $k\rightarrow \infty$.
\end{itemize}
\end{theorem}

In the case of $p=1$, there arises a natural question:   what conditions
should be imposed in Theorem \ref{thapro} on the operator $\left(L_0, B_0\right)$ in order
to assure that  $\left(L_k, B_k\right) \stackrel{s}{\longrightarrow} \left(L_0, B_0\right)$
in the uniform operator topology, in addition to the strong convergence? Note that this is
always true in the case  $1<p<\infty$. The answer to the question is given by the following theorem.

\begin{theorem}\label{thapro2}
Let the assumptions of Theorem \ref{thapro} be satisfied and $p=1$.  Then, there holds
$\big\| \left(L_k, B_k\right)^{-1}-\left(L_0, B_0\right)^{-1}\big\|\rightarrow 0$, as $k\rightarrow \infty$,
if and only if each entry of the matrix-valued function $\Phi_0$ is equal to a regulated function
almost everywhere.
\end{theorem}

Recall that a function on an interval $[a,b]$ is called regulated if it has finite one-sided
limits at every point of the interval (see, for example, \cite{Bourbaki}).

\textbf{Example} Let us consider the Sturm--Liouville differential equation on an interval $(0,1)$ with integral boundary conditions of integer and/or fractional orders
\begin{equation}\label{ex1}
Ly=-y'' + q(t)y=f(t), \quad t\in(0,1),
\end{equation}
\begin{equation}\label{ex2}
By=\sum _{j=1}^{N} \int_{a}^{b}\beta_{j}(t)\big(\textmd{D}_{a+}^{(l_j)}y\big)(t)
{\rm d}t=c  \in \mathbb{C}^2.
\end{equation}
Here, the functions $\{q, f\}\subset W_p^1(0,1)$, the matrix functions $\beta_{j}(\cdot)\in \big(L^\infty(0,1)\big)^{1 \times 2}$, $\big(\textmd{D}_{a+}^{(l_j)}y\big)$ are Caputo derivatives of fractional/integer orders
$0 \leq l_1<l_2<\ldots < l_N< 3-1/p$, $1\leq p <\infty$ (see, for example, \cite{Caputo}).
If the characteristic matrix $M(L,B) \in \mathbb{C}^{2 \times 2}$ of a inhomogeneous boundary-value problem \eqref{ex1}, \eqref{ex2} is non-degenerate, then according to \cite[Corollary 2.4]{NovAM2023} this boundary-value problem has a unique solution $y(\cdot) \in W^3_p(0,1)$. Based on Theorem \ref{thapro}, there exists a sequence of boundary-value problems of the form \eqref{ex1} with polynomial functions $q_k$ and $f_k$, and multipoint boundary conditions of a special form
\begin{equation*} B_ky=\alpha_{0,k}y(0)+\alpha_{1,k}y'(0)+\sum_{j=0}^{N(k)}\beta_{j,k}y''(t_{j,k})=c_k,  \end{equation*}
where $\{t_{j,k}\}$ are some rational points from $[0,1]$, such that,
$\big\|y-y_k\big\|_{3,p} \rightarrow 0$, as $k \rightarrow \infty$,
if $\big\|q-q_k\big\|_{1,p} \rightarrow 0$ and $c_k \rightarrow c$. Moreover, if $p>1$, then
\begin{equation*}
\big\|(L,B)^{-1}-(L_k,B_k)^{-1}\big\| \rightarrow 0, \quad k \rightarrow \infty.
\end{equation*}

\begin{acknowledgement}
The work of the first named author was funded by Postdoctoral Fellowship EU-MSCA4Ukraine (number: 1244691,
WBS-number: 4100609). This project has received funding through the MSCA4Ukraine project, which is funded by
the European Union. Views and opinions expressed are however those of the author only and do not necessarily
reflect those of the European Union, the European Research Executive Agency or the MSCA4Ukraine Consortium.

Neither the European Union nor the European Research  Executive Agency, nor the MSCA4 Ukraine Consortium as
whole nor any individual member institution of the MSCA4Ukraine Consortium can be held responsible for them.

The second named author would like to thank the Isaac Newton Institute for Mathematical Sciences, Cambridge, for support and hospitality during the "Solidarity Program" where work on this paper was undertaken. This work was supported by "EPSRC grant no EP/R014604/1". He also wishes to thank the Department of Mathematics, King's College London, for their hospitality, and to the Ministry of Education and Sciences of Ukraine for support under the grant 0126U000898.
\end{acknowledgement}

\end{document}